\documentclass[10pt,a4paper,oneside]{scrartcl}

\usepackage[utf8]{inputenc}
\usepackage{amsmath}
\usepackage{amsfonts}
\usepackage{amssymb}
\usepackage{epsfig}
\usepackage{graphics} 
\usepackage{hyperref}

\newtheorem{defi}{Definition}[section]
\newtheorem{thm}[defi]{Theorem}

\numberwithin{equation}{section}
\newenvironment{mproof}{\paragraph{Proof.}}{\hfill$\blacksquare$}

\renewcommand{\vec}{\boldsymbol}

\newcommand{\ud}{\,\mathrm{d}}
\newcommand{\R}{\mathbb{R}}
\newcommand{\N}{\mathbb{N}}

\newcommand{\veps}{\vec \varepsilon}
\newcommand{\la}{\langle}
\newcommand{\ra}{\rangle}

\begin{document}

\title{\Large Iterative Coupling of Mixed and Discontinuous\\ Galerkin Methods for 
Poroelasticity}

\date{}


\author{Markus Bause\\ {\small Helmut Schmidt University, Faculty of 
Mechanical Engineering,}\\ {\small Holstenhofweg 85, 220433 Hamburg, 
Germany}}

\pagestyle{myheadings}
\markboth{M.\ Bause}{Mixed and Discontinuous Galerkin Methods for Poroelasticity}

\maketitle

\vspace*{-0.7cm}
\begin{abstract}
\textbf{Abstract.} We analyze an iterative coupling of mixed and discontinuous Galerkin 
methods 
for numerical modelling of coupled flow and mechanical deformation in porous media. The 
iteration is based on an optimized fixed-stress split along with a discontinuous 
variational time discretization. For the spatial discretization of the subproblem of flow 
mixed finite element techniques are applied. The discretization of the subproblem of 
mechanical 
deformation uses discontinuous Galerkin methods. They have shown their ability to 
eliminate locking that sometimes arises in numerical algorithms for poroelasticity and 
causes nonphysical pressure oscillations.
\end{abstract}

\section{Introduction and mathematical model}
\label{Sec:Intro}

We consider the quasi-static Biot system of flow in deformable porous media,
\begin{align}
\label{Eq:B_1}
&- \nabla \cdot \left(2\mu \vec \varepsilon (\vec u) + \lambda \nabla \cdot \vec u \vec I 
- b\,p\vec I\right) = \rho_b \vec g\,,\\[0ex]
\label{Eq:B_2}
&\partial_t\Big(\frac{1}{M}p + \nabla \cdot (b \vec u)\Big) + \nabla \cdot 
 \vec q  = f \,, \quad \vec q = - \vec K \nabla p \,,
\end{align}
in $\Omega\times I$ for a bounded Lipschitz domain $\Omega \subset \R^{d}$, with 
$d = 2,3$, 
and $I=(0,T]$. For simplicity, it is supplemented by homogeneous initial and Dirichlet 
boundary conditions for $p$ and $\vec u$. We 
denote by $\vec u$ the unknown displacement field, $\vec \varepsilon(\vec u) = 
(\nabla \vec u +(\nabla \vec u)^\top)/2$ the linearized strain tensor, $p$ the unknown 
fluid pressure, $\mu$ and $\lambda$ the Lam\'{e} constants, $b$ Biot's coefficient, 
$\rho_b$ the bulk density, $M$ Biot's modulus, $\vec q$ Darcy's velocity and $\vec K$ 
the permeability field. Further, $\vec g$ denotes gravity or, in general, some 
body force and $f$ is a volumetric source. We assume that $\vec g(0) = \vec 0$. The 
quantities $\mu$, $\lambda$, $b$, $\rho_b$ and $M$ are positive constants. The matrix 
$\vec K$ is supposed to be symmetric and uniformly positive definite. Well-posedness of 
\eqref{Eq:B_1}, \eqref{Eq:B_2} has been shown in the literature. 

Recently, iterative coupling schemes for solving the Biot system 
\eqref{Eq:B_1}, \eqref{Eq:B_2} have attracted researchers' interest and have shown their 
efficiency; cf., e.g., \cite{BRK17,BBNKR17,MW13} and the references therein. Iterative coupling  offers the appreciable advantage over the fully coupled 
method that existing and highly developed discretizations and algebraic solver 
technologies can be reused. Here, we use an ''optimized fixed-stress split'' type 
iterative method; cf.\ \cite{BRK17,BBNKR17,MW13}. For the time discretization  of the arising subproblems of flow and mechanical deformation a discontinuous 
in time variational approach (cf.\ \cite{BRK17}) is applied. 

For the approximation of the pressure and flux unknown $p$ and $\vec q$ within the iterative coupling 
approach we use a mixed finite element method. The displacement variable $\vec u$ is 
discretized by a discontinuous Galerkin method. Thereby the results herein represent the 
natural extension of the work of the author in \cite{BRK17} and further work in the 
literature, where the displacement field $\vec u$ was approximated by a continuous 
Galerkin method. The motivation for using a discontinuous Galerkin scheme for the 
discretization of the displacement $\vec u$ comes from combating the locking phenomenon, 
that sometimes arises in numerical algorithms for poroelasticity and manifests as spurious 
nonphysical pressure oscillations. In poroelasticity the locking-dominant parameter is the 
specific storage coefficient $c_0= 1/M$ in Eq.\ \eqref{Eq:B_2}. Locking primarily arises 
if $c_0 = 0$ and, usually, does not appear for $c_0\neq 0$.  For a more extensive discussion of 
locking in poroelasticity we refer to the literature; cf., e.g., 
\cite{PW09,L17}.  

In \cite{PW08} a coupling of mixed and discontinuous Galerkin finite methods is studied 
within a monolithic approach and an error analysis for the semi-discretization is space 
is 
given. By numerical experiments it is shown for the problem of Barry and Mercer, that a  
discontinuous Galerkin discretization of the displacement $\vec u$ is capable of 
eliminating spurious pressure oscillations related to locking arising in continuous 
discretizations of $\vec u$. Therefore, it seems worth to study the combined mixed 
and discontinuous Galerkin approach also within a fixed-stress split iterative method 
which is done here.

We use standard notation. In particular, we put $W=L^2(\Omega)$ and $\vec V = \vec 
H(\text{div};\Omega)$ and denote by $\langle \cdot, \cdot\rangle$ the inner product of 
$W$.

\section{Iterative coupling and space-time discretization}
\label{Sec:Disc}

We consider solving the system \eqref{Eq:B_1}, \eqref{Eq:B_2} by the 
following iteration scheme.

\medskip
\textbf{Subproblem of flow:} \textit{Let $ f^{\,k} := f - b \, \nabla \cdot \partial_t 
\vec u^k 
+ L\, \partial_t p^k \in 
L^2(I;W)$ be given. Find $p^{k+1}\in H^1(I;W)$, $\vec q^{k+1}\in L^2(I;\vec V)$ such 
that $p^{k+1}(0)=0$ and}
\begin{align}
\label{Eq:WPF1}
& \mbox{$\left(\frac{1}{M}+ L \right) \int_I \langle \partial_t p^{k+1}, w\rangle \ud t + 
\int_I \langle \nabla \cdot \vec q^{k+1}, w\rangle \ud t$} =  
\mbox{$\int_I \langle  f^{\,k}, w \rangle \ud t$}\,, \\[1ex]
\label{Eq:WPF2}
& \mbox{$\int_I \langle \vec K^{-1} \vec q^{k+1}, \vec v\rangle \ud t - \int_I \langle 
p^{k+1}, \nabla \cdot \vec v \rangle \ud t$}  = 0 
\end{align}
\textit{for all $w\in L^2(I;W)$ and $\vec v\in L^2(I;\vec V)$.} 

\smallskip
\textbf{Subproblem of mechanical deformation:} \textit{Let $p^{k+1}\in 
H^1(I;W)$ be given. Find $\vec u^{k+1} \in H^1(I;\vec H^1(\Omega))\cap L^2(I;\vec 
H^1_0(\Omega))$ such that $\vec 
u(0)=\vec 0$ and}
\begin{equation}
\label{Eq:WPM}
\begin{aligned}
&\mbox{$\int_I 2 \mu \langle \vec \varepsilon(\vec u^{k+1}), \vec \varepsilon(\vec 
z)\rangle \ud 
t + \int_I \lambda \langle  \nabla \cdot \vec u^{k+1}, \nabla \cdot \vec z \rangle \ud 
t$}\\[1ex]
&\qquad \mbox{$ = \rho_b \int_I \langle \vec g , \vec z \rangle \ud t +  b \int_I \langle 
p^{k+1}, 
\nabla \cdot \vec z\rangle \ud t$}
\end{aligned}
\end{equation}
\textit{for all $\vec z \in L^2(I;\vec H^1_0(\Omega))$.}

\smallskip
In this scheme the artificial quantity $L$ is a numerical parameter that was 
firstly introduced in \cite{MW13} to accelerate the iteration process. The 
convergence of the iteration \eqref{Eq:WPF1}--\eqref{Eq:WPM} is ensured for all $L\geq 
b^2/(2\lambda)$; cf.\ \cite[Thm.\ 2.1]{BRK17}.

For the discretization we decompose the time interval $(0,T]$ into $N$ subintervals 
$I_n=(t_{n-1},t_n]$, where $n\in \{1,\ldots ,N\}$ and $0=t_0<t_1< \cdots < t_{N-1} < t_N 
= T$ and $\tau = \max_{n=1,\ldots N} (t_n-t_{n-1})$. We denote by $\mathcal{T}_h=\{K\}$ a 
finite element decomposition of mesh size $h$ of the polyhedral domain 
$\overline{\Omega}$ into closed subsets $K$, quadrilaterals in two dimensions and 
hexahedrals in three dimensions. Further, $\mathcal E_{h,}^{\text{int}}$ is the set of 
all interior edges (faces for $d=3$). To each interior edge (or face) $e\in \mathcal 
E_{h}^{\text{int}}$ we associate a fixed unit normal vector $\vec \nu^e$.

For the spatial discretization of \eqref{Eq:WPF1}, 
\eqref{Eq:WPF2} we use a mixed finite element approach. We choose the class of 
Raviart–Thomas elements for the two-dimensional case and 
the class of Raviart--Thomas--N\'ed\'elec elements in three space
dimensions, where for $s\geq 0$ the space $W_h^s \subset W$ with $
W_{h}^s=\{ w_h\in L^2(\Omega) \mid w_h{}_{|_{K}}\circ T_K \in \mathbb{Q}_s\}$ and $\vec 
V_h^s \subset \vec V$ denote the corresponding inf-sup stable pair of 
finite element spaces; cf.\ \cite{BRK17}.  Here, $\mathbb Q_s$ is the space of 
polynomials that are of degree less than or equal to 
$s$ in each variable and $T_K$ is a suitable invertible 
mapping of the reference cube $\widehat K$ to $K$ of $\mathcal T_h$.  For the 
spatial discretization of \eqref{Eq:WPM} we discretize we use a discontinuous Galerkin 
method with the space 
\[
 \vec H_{h}^l=\{\vec z_h\in L^2(\Omega) \mid \vec 
z_h{}_{|_{K}}\circ T_K \in \mathbb{Q}_l^d\,, \; \vec z_h{}_{|\partial \Omega} = \vec 
0\}\,.
\] 
The fully discrete space-time finite element spaces are then defined by
\begin{align}
\label{Def:FE4}
{\mathcal{W}}_{\tau,h}^{r,s} & = \{w_{\tau,h}\in L^2( 
I;W)\mid 
w_{\tau,h}{}_{|I_n}\in \mathcal P_r(I_n;W_h^s)\,, \; w_{\tau,h}(0)\in W_h^s\} \,,\\[0ex]
\label{Def:FE5}
\vec{{\mathcal{V}}}_{\tau,h}^{r,s} & = \{\vec v_{\tau,h}\in L^2( 
I;\vec V)\mid 
\vec v_{\tau,h}{}_{|I_n}\in \mathcal P_r(I_n;\vec V_{h}^s)\,,\; \vec v_{\tau,h}(0)\in 
\vec V_h^s\} \,,\\[0ex]
\label{Def:FE6}
\vec{{\mathcal{Z}}}_{\tau,h}^{r,l} & = \{\vec z_{\tau,h}\in 
L^2(I;\vec{L}^2(\Omega))\mid \vec z_{\tau,h}{}_{|I_n}\in \mathcal P_r(I_n;\vec 
H_{h}^l)\,,\; \vec z_{\tau,h}(0)\in \vec H_h^l\}\,,
\end{align}
where $\mathcal P_r(I_n;X)$ is the space of all polynomials in time up to degree 
$r\geq 0$ on $I_n$ with values in $X$. We choose $l=s+1$ to equilibrate the convergence 
rates of the spatial discretization for the three unknowns $p,\vec q$ and 
$\vec u$; cf.\ \cite[p.\ 426, Thm.\ 1]{PW08}. For short, we put $W_h= W_h^s$, $\vec V_h 
= \vec V_{h}^s$ and $\vec H_h = \vec H_{h}^{s+1}$. 

On each time interval $ I_n$ we expand the discrete functions $p^{k}_{\tau,h} \in 
{\mathcal{W}}_{\tau,h}^{r,s}$, $\vec q^{k}_{\tau,h}\in 
\vec{{\mathcal{V}}}_{\tau,h}^{r,s}$ 
and $\vec u^{k}_{\tau,h} \in \vec{{\mathcal{Z}}}_{\tau,h}^{r,s+1}$ in time in 
terms of Lagrangian basis functions $\varphi_{n,j}\in \mathcal P_r(I_n;\R)$ with respect 
to $r+1$ nodal points 
$t_{n,j}\in I_n$, for $j=0,\ldots,r$, 
\begin{equation}
\label{Eq:dGRepSolBasis}
p_{\tau,h}^{k}{}_{|I_n} (t) = \sum_{j=0}^r P^{j,k}_{n,h} 
\varphi_{n,j}(t)\,, \quad \vec z_{\tau,h}^{k}{}_{|I_n} (t) = \sum_{j=0}^r \vec 
Z^{j,k}_{n,h} \varphi_{n,j}(t)\,,
\end{equation}
for $t\in I_n$ and $\vec z \in \{\vec q, \vec u\}$ with coefficient functions $P^{j,k}_{n,h} \in W_h$, $\vec Q^{j,k}_{n,h} \in \vec V_h$ for $ \vec z = \vec q$ and $\vec U^{j,k}_{n,h} \in 
\vec H_h$ for $ \vec z = \vec u$ and $j=0,\ldots ,r$. The nodal points $t_{n,j}$ are 
chosen as the quadrature points of the ($r$+1)-point Gauss quadrature formula on $I_n$ 
which is exact for polynomials of degree less or equal to $2r+1$. 

Solving the variational problems \eqref{Eq:WPF1}, \eqref{Eq:WPF2} and \eqref{Eq:WPM} in 
the fully discrete function spaces \eqref{Def:FE4}--\eqref{Def:FE6} and using a 
discontinuous test basis in time with support on $I_n$ then leads us to the following 
fully discrete iteration scheme, referred to as the MFEM($s$)dG($s$+1)--dG($r$) splitting 
scheme. 

\medskip
\textbf{Fully discrete subproblem of flow:} {\em Let $n\in 
\{1,\ldots ,N\}$. Find coefficient functions $P_{n,h}^{i, k+1}\in W_h$ and $\vec 
Q_{n,h}^{i, k+1}\in \vec V_h$ for $i=0,\ldots, r$ such that}
\begin{align}
\nonumber
&\dfrac{1}{M} \hspace*{-0.3ex} \sum_{j =0}^r \hspace*{-0.3ex}  \alpha_{ij} \langle 
P_{n,h}^{j, k+1},  w_h \rangle 
+ 
L 
\sum_{j =0}^r \hspace*{-0.3ex}  \alpha_{ij} \langle P_{n,h}^{j, k+1} - P_{n,h}^{j, k},   
w_h 
\rangle + 
\tau_n \beta_{ii} \langle \nabla \cdot \vec Q_{n,h}^{i, k+1} , w_h \rangle \\[-1ex]
\label{Eq:D_AlgebProb1}
&  = 
\tau_n \beta_{ii}\langle f(t_{n,i}), w_h\rangle -b \sum_{j =0}^r 
\alpha_{ij} \sum_{K\in \mathcal T_h} \langle \nabla 
\cdot \vec U_{n,h}^{j, k} , w_h \rangle_K 
\\[-1ex] \nonumber
&+ \gamma_i \, \frac{1}{M}\langle p^\infty_{\tau,h}(t_{n-1}^-),w_h\rangle + \gamma_i 
\, b \sum_{K\in \mathcal T_h}\langle  \nabla \cdot \vec 
u^\infty_{\tau,h}(t_{n-1}^-),w_h\rangle_K\,,\\[0ex]
\label{Eq:D_AlgebProb2}
&\langle \vec K^{-1} \vec Q_{n,h}^{i, k+1}, \vec v_h \rangle - \langle P_{n,h}^{i, k+1}, 
\nabla \cdot \vec v_h \rangle   = 0
\end{align}
\emph{for all $w_h \in W_h$, $\vec v_h \in \vec V_h$ and $i= 0,\ldots, r$, where 
$p_{\tau,h}(t_{n-1}^-)=0$ for $n=1$ and $p^\infty_{\tau,h}(t_{n-1}^-) = 
\lim_{k\rightarrow \infty}p^k_{\tau,h}(t_{n-1}^-)$ for $n>1$, and similarly 
for $\vec u^\infty_{\tau,h}(t_{n-1}^-)$.}

\medskip
\textbf{Fully discrete subproblem of mechanical deformation:} {\em Let $n\in 
\{1,\ldots ,N\}$. Find coefficient functions $\vec U_{n,h}^{i, k+1} \in \vec H_h$ for 
$i=0,\ldots, r$ such that}
\begin{equation}
\label{Eq:D_AlgebProb3}
\begin{aligned}
& \sum_{K\in \mathcal T_h} \langle \vec \sigma(\vec U_{n,h}^{i, k+1}),\vec 
\varepsilon(\vec z_h )\rangle_K + J_\delta (\vec U_{n,h}^{i, k+1},\vec z_h)- J_d (\vec 
\sigma(\vec U_{n,h}^{i, 
k+1}),\vec z_h) \\[0ex]
& = \sum_{K\in \mathcal T_h} \ b \langle P_{n,h}^{i, k+1}, \nabla \cdot \vec z_h 
\rangle_K - b \, J_p(P_{n,h}^{i, k+1},\vec z_h)
\end{aligned}
\end{equation}
\emph{for all $\vec z_h \in \vec H_h$ and $i= 0,\ldots, r$ with effective stress 
$\vec \sigma(\vec u)= 2\mu \vec \varepsilon (\vec u) + \lambda \nabla \cdot \vec u \, 
\vec I$.}

\smallskip
In Eq.\ \eqref{Eq:D_AlgebProb3} we use the notation
\begin{align*}
J_\delta(\vec y_h,\vec z_h) & = \sum_{e \in \mathcal 
E_{h}^{\text{int}}} \frac{\delta_e}{|e|^\beta} \langle [\vec y_h],[\vec z_h]\rangle_e \,, 
\quad J_p(w_h,\vec z_h)= \sum_{e \in \mathcal 
E_{h}^{\text{int}}} \langle \{w_h\} \, \vec \nu^e, [\vec z_h]\rangle_e\,,\\[0ex]
J_{d}(\vec y_h,\vec z_h) & = \sum_{e \in \mathcal 
E_{h}^{\text{int}}} \big( \langle \{\vec \sigma(\vec y_h)\vec \nu^e\}, [\vec 
z_h]\rangle_e +  \langle \{\vec \sigma(\vec z_h)\vec 
\nu^e\}, [\vec y_h]\rangle_e\big)\,.
\end{align*}
We denote by $\langle \cdot , \cdot \rangle_K$ and $\langle \cdot , \cdot \rangle_e$ the 
$L^2$ inner products on $K$ and $e$, respectively. As usual, we let $\{w\}= 
((w_{|K})_{|e}+(w_{|K'})_{|e})/2$ for two adjacent elements $K$ and $K'$ with common edge 
(or face) $e$ and, similarly, $[w]= (w_{|K})_{|e}-(w_{|K'})_{|e}$. The penalty term 
$J_\delta$ contains the numerical parameter $\delta_e$ that takes a 
constant value at each edge (or face) $e\in \mathcal E_h^{\text{int}}$ with Lebesgue 
measure $|e|$. The power $\beta$ is a positive number that depends on the dimension 
$d$. In \cite{PW08}, the choice $\beta=(d-1)^{-1}$ is proposed for the fully coupled 
semi-discretization of \eqref{Eq:B_1}, \eqref{Eq:B_2}. The coefficients $\alpha_{ij}$, 
$\beta_{ii}$ and $\gamma_i$ 
are defined by $\alpha_{ij} = \int_{I_n} \varphi_{n,j}'(t) \cdot 
\varphi_{n,i}(t)\ud t+ \gamma_i\cdot\gamma_j$, $\beta_{ii} = \int_{I_n} 
\varphi_{n,i} (t) \cdot \varphi_{n,i}(t)\ud t$ and $\gamma_{i} = 
\varphi_{n,i}(t_{n-1}^+)$ for $i,j=0,\ldots,r$. For a detailed derivation of the 
Eqs.\ \eqref{Eq:D_AlgebProb1}--\eqref{Eq:D_AlgebProb3} we refer to \cite{BRK17} 
for the space-time issue and to \cite{PW08} for the discontinuous Galerkin approximation 
in space of the displacement field $\vec u$. Eq.\ \eqref{Eq:D_AlgebProb3} is referred to 
as the symmetric interior penalty (SIP) discontinuous Galerkin method. 
Finally, we note that an additional penalty term involving $\partial_t \vec u$ is 
proposed for the semi-discrete fully coupled approach in \cite{PW08}.

\section{Convergence of the iteration scheme}
\label{Sec:Conv}

Here we prove the convergence of the iteration scheme \eqref{Eq:D_AlgebProb1}, \eqref{Eq:D_AlgebProb2} and \eqref{Eq:D_AlgebProb3}. For this, let $\{p_{\tau,h},\vec q_{\tau,h}\} \in 
\mathcal{W}_{\tau,h}^{r,s} \times \vec{\mathcal{V}}_{\tau,h}^{r,s}$, $\vec u_{\tau,h} 
\in 
\vec{\mathcal{Z}}_{\tau,h}^{r,s+1}$ denote the fully discrete 
MFEM($s$)dG($s$+1)--dG($r$) approximation of 
\eqref{Eq:B_1}, \eqref{Eq:B_2}, formally given by passing to the limit $k\rightarrow 
\infty$ in the scheme \eqref{Eq:D_AlgebProb1}--\eqref{Eq:D_AlgebProb3}. Analogously to 
\eqref{Eq:dGRepSolBasis}, let $\{p_{\tau,h} ,\vec q_{\tau,h},\vec u_{\tau,h}\}$ 
on $I_n$ be represented by coefficients $P_{n, h}^{j}\in W_h$, $\vec Q_{n,h}^j\in \vec 
V_h$ and $\vec U_{n,h}^j\in \vec H_h$ for $j=0,\ldots, r$. For $n\in \{1,\ldots ,N\}$ and 
$t\in I_n$ we put  
\begin{equation*}
\mbox{$E_p^{j, k} = P_{n, h}^{j, k} -  P_{n, h}^{j}\,, \quad
e_p^{k}(t) = \sum_{j=0}^r E_p^{j, k}\varphi_{n,j}(t)\,, \quad S_p^{i, k} = \sum_{j 
=0}^r \alpha_{ij} E_p^{j, k}$}\,.
\end{equation*}
The quantities $\vec E_{\vec q}^{j, k}$, $\vec e_{\vec q}^{k}$, $\vec E_{\vec u}^{j, 
k}$, $\vec e_{\vec u}^{k}$ and $\vec S_{\vec q}^{i, k}$, $\vec 
S_{\vec u}^{i, k} $ are defined analogously. 

\begin{thm} 
\label{Thm:ConvDisc}
Let $\{p_{\tau,h},\vec q_{\tau,h}\} \in 
\mathcal{W}_{\tau,h}^{r,s} \times \vec{\mathcal{V}}_{\tau,h}^{r,s}$, $\vec u_{\tau,h} 
\in 
\vec{\mathcal{Z}}_{\tau,h}^{r,s+1}$ denote the fully discrete 
space-time MFEM($s$)dG($s$+1)--dG($r$) approximation of 
\eqref{Eq:B_1}, \eqref{Eq:B_2}. Let $\{p^k_{\tau,h} ,\vec 
q^k_{\tau,h}, \vec u^k_{\tau,h} \}$ be defined  by \eqref{Eq:dGRepSolBasis} with 
coefficients being given by \eqref{Eq:D_AlgebProb1}, \eqref{Eq:D_AlgebProb2} and 
\eqref{Eq:D_AlgebProb3}. Then, if the parameter $L > 0$ and penalty function $\delta$ in 
$J_\delta$ of \eqref{Eq:D_AlgebProb3} are chosen sufficiently large, the sequence 
$\{S_{p}^{i,k}\}_k$, for $i=0,\ldots,r$, converges in $W_h$. This implies the convergence 
of $\{p^k_{\tau,h}(t_n^\pm),\vec q^k_{\tau,h}(t_n^\pm),\vec u^k_{\tau,h}(t_n^\pm)\}$ to 
$ \{p_{\tau,h}(t_n^\pm),\vec q_{\tau,h}(t_n^\pm),\vec u_{\tau,h}(t_n^\pm)\}$ in 
$W_h\times \vec V_h \times \vec H_h$ for $k\rightarrow \infty$ and $n=1,\ldots, N$, as 
well as of  $p^k_{\tau,h}$, $\vec q^k_{\tau,h}$ and $\vec u^k_{\tau,h}$ in  $L^2(I_n;W)$ 
and $L^2(I_n;\vec L^2(\Omega))$, respectively. 
\end{thm}

\begin{mproof}
We split the proof into several steps.

\smallskip
\noindent \textbf{1.\ Step (Error equations).} By substracting Eqs.\ 
\eqref{Eq:D_AlgebProb1}--\eqref{Eq:D_AlgebProb3} from  the fully discrete 
monolithic space-time approximation MFEM($s$)dG($s$+1)--dG($r$) of the Biot system 
\eqref{Eq:B_1}, \eqref{Eq:B_2} we obtain for $i=0, \ldots, r$ that
\begin{align}
& \dfrac{1}{M} \sum_{j =0}^r \alpha_{ij} \la E_p^{j, k+1}, w_h \ra + L \sum_{j =0}^r 
\alpha_{ij} \la E_p^{j, k+1} - E_p^{j, k}, w_h \ra \qquad \qquad \qquad \nonumber 
\\[-1.5ex]
& \qquad + \tau_n \beta_{ii} \la \nabla \cdot \vec E_{\vec q}^{i, k+1}, w_h \ra = 
-b \sum_{j 
=0}^r \alpha_{ij} \sum_{K\in \mathcal T_h} \la \nabla \cdot \vec E_{\vec u}^{j, k}, w_h 
\ra_K \,,
\label{proof_eq_7}\\[0ex]
& \la \vec K^{-1} \vec E_{\vec q}^{i, k+1}, \vec v_h \ra - \la E_p^{i, k+1}, \nabla \cdot 
\vec v_h \ra = 0\,,
\label{proof_eq_8}\\[0.5ex]
& 
\begin{aligned}
&\sum_{K\in \mathcal T_h}  \langle \vec \sigma(\vec E_{\vec u}^{i, k+1}),\vec 
\varepsilon(\vec z_h )\rangle_K + J_\delta(\vec E_{\vec u}^{i, k+1},\vec z_h)- 
J_d(\vec E_{\vec u}^{i, k+1},\vec z_h)\\[-1ex]
& \qquad  = \sum_{K\in \mathcal T_h}  b \langle E_{p}^{i, k+1}, \nabla \cdot \vec 
z_h \rangle_k - b \, J_p(E_{p}^{i, k+1},\vec z_h)
\end{aligned}
\label{proof_eq_9}
\end{align}
for all $w_h \in W_h$, $\vec v_h \in \vec V_h$, $\vec z_h \in \vec H_h$.

\smallskip
\noindent\textbf{2.\ Step (Choice of test function in Eq.\ 
{\eqref{proof_eq_7}}).} 
We test Eq.~\eqref{proof_eq_7} with $w_h = \sum_{j =0}^r \alpha_{ij} E_p^{j, k+1}$. By 
some calculations following \cite[p.\ 756]{BRK17} we get that  
\begin{equation}
\label{proof_eq_11}
\begin{aligned}
& \Big(\dfrac{1}{M \beta_{ii}}  + \dfrac{L}{ 2 \beta_{ii}}\Big)  \| S_p^{i, k+1} \|^2  
+ 
\dfrac{L}{ 2 \beta_{ii}}  \| S_p^{i, k+1} - S_p^{i, k} \|^2 - \dfrac{L}{ 2 \beta_{ii}} 
\| 
S_p^{i, k} \|^2 \\[0ex]
& + \tau_n \la \nabla \cdot \vec E_{\vec q}^{i, k+1},  S_p^{i, k+1} \ra = - 
\dfrac{b}{\beta_{ii}} \sum_{K\in \mathcal T_h}  \la \nabla \cdot \vec S_{\vec u}^{i, k}, 
S_p^{i, k+1} \ra_K\,.
\end{aligned}
\end{equation}

\smallskip
\noindent\textbf{3.\ Step (Summation of Eq.\ \eqref{proof_eq_8} and choice of test 
function).} Changing the index $i$ in Eq.\  \eqref{proof_eq_8} to $j$, multiplying this 
equation with $\alpha_{ij}$, summing up from $j=0$ to $r$ and testing with $\vec v_h 
=\tau_n \vec E_{\vec q}^{i, k+1} \in \vec V_h$ we have that
\begin{equation}\label{proof_eq_14}
\tau_n \la \vec K^{-1} \vec S_{\vec q}^{i, k+1}, \vec E_{\vec q}^{i, k+1} \ra - \tau_n 
\la S_p^{i, k+1}, \nabla \cdot \vec E_{\vec q}^{i, k+1}\ra = 0\,.
\end{equation}
Adding Eq.\ 
\eqref{proof_eq_14} to Eq.\ \eqref{proof_eq_11} implies that
\begin{equation}
\label{proof_eq_15}
\begin{aligned}
& \Big(\dfrac{1}{M \beta_{ii}} + \dfrac{L}{ 2 \beta_{ii}}\Big)\| S_p^{i, k+1} \|^2  + 
\dfrac{L}{ 2 \beta_{ii}} \| S_p^{i, k+1} - S_p^{i, k} \|^2 - \dfrac{L}{ 2 \beta_{ii}} \| 
S_p^{i, k} 
\|^2 \\[0ex]
& + \tau_n \la \vec K^{-1} \vec S_{\vec q}^{i, k+1}, \vec E_{\vec q}^{i, k+1}\ra = - 
\dfrac{b}{\beta_{ii}} \sum_{K\in \mathcal T_h}  \la \nabla \cdot \vec S_{\vec u}^{i, k}, 
S_p^{i, k+1} \ra_K \,.
\end{aligned}
\end{equation}

\smallskip
\noindent\textbf{4.\ Step (Summation of Eq.\ \eqref{proof_eq_9} and choice of test 
function).} Changing the index $i$ in Eq.\ 
\eqref{proof_eq_9} to $j$, multiplying the resulting equation with $\alpha_{ij}$, summing 
up from $j=0$ to $r$ and choosing $\vec z_h = \beta_{ii}^{-1} 
\vec S_{\vec u}^{i, k} \in \vec H_h$ we find that
\begin{align}
 \nonumber
&\frac{1}{\beta_{ii}}\sum_{K\in \mathcal T_h}  \langle \vec \sigma(\vec S_{\vec u}^{i, 
k+1}),\vec 
\varepsilon(\vec S_{\vec u}^{i, k}  )\rangle_K + \frac{1}{\beta_{ii}} J_\delta(\vec 
S_{\vec 
u}^{i, k+1} ,\vec S_{\vec u}^{i, k} ) - \frac{1}{\beta_{ii}} J_d(\vec S_{\vec 
u}^{i, k+1} ,\vec S_{\vec u}^{i, k} )\\[-0.5ex]
& = \frac{b}{\beta_{ii}} \sum_{K\in \mathcal T_h}  \langle S_{p}^{i, k+1}, \nabla \cdot 
\vec S_{\vec u}^{i, k}  \rangle_K - \frac{b}{\beta_{ii}} J_p(S_{p}^{i, k+1}, \vec 
S_{\vec u}^{i, k})\,.
\label{proof_eq_13}
\end{align}
Adding Eq.\ \eqref{proof_eq_13} to Eq.\ \eqref{proof_eq_15} leads to
\begin{align}
\nonumber
& \Big(\dfrac{1}{M \beta_{ii}}  + \dfrac{L}{ 2 \beta_{ii}} \Big)\| S_p^{i, k+1} \|^2 
+ \dfrac{L}{ 2 \beta_{ii}} \| S_p^{i, k+1} - S_p^{i, k} \|^2 + \tau_n \la \vec K^{-1} \vec 
S_{\vec q}^{i, k+1}, \vec E_{\vec q}^{i, k+1}\ra\\
\label{proof_eq_17}
& \quad + \frac{1}{\beta_{ii}} J_\delta(\vec S_{\vec u}^{i, k+1} ,\vec S_{\vec u}^{i, k}) + 
\frac{1}{\beta_{ii}}\sum_{K\in \mathcal T_h} \la \vec \sigma(\vec S_{\vec u}^{i, k+1}), 
\vec \veps(\vec S_{\vec u}^{i, k})\ra_K \\[0ex]
\nonumber
& =  \dfrac{L}{ 2 \beta_{ii}} \| S_p^{i, k} \|^2 - 
\frac{b}{\beta_{ii}} J_p(S_{p}^{i, k+1},\vec S_{\vec u}^{i, k}) + \frac{1}{\beta_{ii}}  
J_d(\vec S_{\vec 
u}^{i, k+1} ,\vec S_{\vec u}^{i, k})\,.
\end{align}

\smallskip
\noindent\textbf{5.\ Step (Formation of incremental equation for \eqref{proof_eq_9}, 
summation and choice of test function.)} 
Firstly, we write Eq.\ \eqref{proof_eq_9} for two consecutive iterations, $k$ and 
$k +1$, and substract the resulting equations from each other. Secondly, we change the 
index $i$ in the thus obtained equations to $j$, multiply them with $\alpha_{ij}$ and sum 
up from $j=0$ to $r$ to obtain that 
\begin{align}
\nonumber 
& \sum_{K\in \mathcal T_h} \la \vec \sigma (\vec S_{\vec 
u}^{i, k+1} -\vec S_{\vec u}^{i, k}), \vec \veps(\vec z_h)\ra_K +J_\delta (\vec S_{\vec 
u}^{i, k+1}-\vec S_{\vec u}^{i, k},\vec z_h)= b \sum_{K\in \mathcal T_h}  \la  S_p^{i, 
k+1} \\
\label{proof_eq_17b}
&   - S_p^{i, k}, \nabla \cdot \vec z_h  \ra_K - {b} J_p(S_{p}^{i, k+1}-S_{p}^{i, k},\vec 
z_h) + J_d(\vec S_{\vec u}^{i, k+1}-\vec 
S_{\vec u}^{i, k},\vec z_h)
\end{align}
for all $\vec z_h \in \vec H_h$.
Choosing $ \vec z_h = \vec S_{\vec u}^{i, k+1} - \vec S_{\vec u}^{i, k} \in \vec H_h $ in 
\eqref{proof_eq_17b}, dividing by $ \beta_{ii} > 0$ and  summing up the resulting 
identity from $i = 0 $ to $r$, we find that 
\begin{align}
\nonumber
& \sum_{i=0}^r\frac{1}{\beta_{ii}} \Big(\hspace*{-0.5ex}\sum_{K\in \mathcal T_h} 
\hspace*{-1ex} \la \vec \sigma (\vec S_{\vec 
u}^{i, k+1} \hspace*{-0.3ex}-\vec S_{\vec u}^{i, k}), \vec \veps(\vec S_{\vec 
u}^{i, k+1}\hspace*{-0.3ex} -\vec S_{\vec u}^{i, k})\ra_K + J_\delta (\vec S_{\vec u}^{i, 
k+1}\hspace*{-0.3ex}-\vec S_{\vec u}^{i, k},\vec S_{\vec u}^{i, k+1}
\end{align}

\begin{align}
\label{proof_eq_18}
&  - \vec S_{\vec u}^{i, k})\Big) = \sum_{i=0}^r 
\frac{1}{\beta_{ii}} \Big(\sum_{K\in \mathcal T_h} b \la  S_p^{i, k+1} - 
S_p^{i, k}, \nabla \cdot (\vec S_{\vec 
u}^{i, k+1} -\vec S_{\vec u}^{i, k})  \ra_K \hspace*{-3cm}\mbox{}\\[0ex]
\nonumber
& -  b\, J_p(S_{p}^{i, k+1}-S_{p}^{i, k},\vec 
S_{\vec u}^{i, k+1} - \vec S_{\vec u}^{i, k})
+ J_d(\vec S_{\vec u}^{i,k+1}-\vec S_{\vec u}^{i,k},\vec 
S_{\vec u}^{i, k+1} - \vec S_{\vec u}^{i, k} )\Big) \,.
\end{align}
Further, from Eq.\ \eqref{proof_eq_17b} with $ \vec z_h = \vec 
S_{\vec u}^{i, k+1} - \vec S_{\vec u}^{i, k}$ we get by means of the inequalities of 
 Cauchy--Schwarz and Cauchy--Young that
\begin{align}
\nonumber
& \sum_{K\in \mathcal T_h} \la \vec \sigma (\vec S_{\vec 
u}^{i, k+1} \hspace*{-0.3ex}-\vec S_{\vec u}^{i, k}), \vec \veps(\vec S_{\vec 
u}^{i, k+1}\hspace*{-0.3ex} -\vec S_{\vec u}^{i, k})\ra_K +J_\delta (\vec S_{\vec u}^{i, 
k+1}-\vec S_{\vec u}^{i, k},\vec 
S_{\vec u}^{i, k+1} \\[-1ex]
\label{proof_eq_19}
&  - \vec S_{\vec u}^{i, k}) \leq \frac{b^2}{2\lambda}  \|  S_p^{i, k+1} 
- S_p^{i, k}\|^2 + \sum_{K\in \mathcal T_h} \frac{\lambda}{2} \|\nabla 
\cdot (\vec S_{\vec u}^{i, k+1} - \vec S_{\vec u}^{i, k})\|_K^2\\[0ex]
\nonumber
& - {b} \, J_p(S_{p}^{i, k+1}-S_{p}^{i, k},\vec S_{\vec u}^{i, k+1} - \vec S_{\vec u}^{i, 
k}) + J_d(\vec S_{\vec u}^{i, k+1}- \vec S_{\vec u}^{i, k},\vec 
S_{\vec u}^{i, k+1} - \vec S_{\vec u}^{i, k})\,.
\end{align}
The terms $J_p$ and $J_d$ can be bounded by means of (cf.\ \cite[p.\ 429, p.\ 431]{PW08}) 
\begin{equation}
\label{Eq:BoundJ}
\hspace*{-3cm}
\begin{aligned}
|J_p(w_h,\vec z_h) | & \leq  \frac{c}{\delta_{\min}} \|w_h\|^2 + \frac{1}{R} J_\delta 
(\vec z_h,\vec z_h)\,,\\
\Big| \sum_{e \in \mathcal 
E_{h}^{\text{int}}} \langle \{\vec \sigma(\vec y_h)\vec \nu^e\}, [\vec 
z_h]\rangle_e\Big| & \leq \frac{c}{\delta_{\min}} \sum_{e\in \mathcal 
E_h^{\text{int}}} \langle \vec \sigma (\vec y_h), \vec \varepsilon (\vec y_h)\rangle_K + 
\frac{1}{R} J_\delta (\vec z_h,\vec z_h)
\end{aligned}
\hspace*{-3cm}\mbox{}
\end{equation}
for $R\in \N$, $R>1$ and some constant $c>0$, such that they can be absorbed by the 
left-hand side, if the penalty parameter $\delta_{\min}$ is chosen sufficiently large.  

\smallskip
\textbf{6.\ Step (Summation of Eq.\ \eqref{proof_eq_17} and combination 
with Eq.\ \eqref{proof_eq_18}).} 
Using in \eqref{proof_eq_17} that $4\la x, y \ra = \| x + y \|^2 - \| 
x - y \|^2$,  summing up the resulting equation over $i$ and  
using \eqref{proof_eq_18} together with \eqref{Eq:BoundJ}) we get that
\begin{align}
\nonumber
& \sum_{i=0}^r \Big\{ \Big(\dfrac{1}{M \beta_{ii}} + \dfrac{L}{ 
2 \beta_{ii}}\Big) \| S_p^{i, k+1} \|^2  + \dfrac{L}{ 2 \beta_{ii}} \| S_p^{i, 
k+1} - S_p^{i, k} \|^2  \\[-1ex]
\nonumber 
& + \tau_n \la \vec K^{-1} \vec S_{\vec q}^{i, k+1}, \vec E_{\vec 
q}^{i, k+1}\ra\Big\} + \sum_{i=0}^r \frac{1}{4\beta_{ii}} J_\delta(\vec 
S_{\vec u}^{i, k+1} + \vec S_{\vec u}^{i, k},\vec 
S_{\vec u}^{i, k+1} + \vec S_{\vec u}^{i, k} )\\[-1ex]
\nonumber
&  + \sum_{i=0}^r\frac{1}{4\beta_{ii}} \sum_{K\in \mathcal T_h} \la \vec \sigma (\vec S_{\vec 
u}^{i, k+1} +\vec S_{\vec u}^{i, k}), \vec \veps(\vec S_{\vec 
u}^{i, k+1} +\vec S_{\vec u}^{i, k})\ra_K\leq   \sum_{i=0}^r  
\dfrac{L}{ 2 \beta_{ii}} \| S_p^{i, k} \|^2  \\[-1ex]
\nonumber
& + \sum_{i=0}^r 
\dfrac{1}{4 \beta_{ii}} \sum_{K\in \mathcal T_h} b \la S_p^{i, k+1} - S_p^{i, k}, 
\nabla \cdot (\vec S_{\vec u}^{i, k+1} - \vec S_{\vec u}^{i, k}) \ra_K+ c \sum_{i=0}^r 
\frac{b^2}{\beta_{ii}}\|S_p^{i, k+1}\\[-1ex]
\label{proof_eq_22}
& - S_p^{i, k}\|^2 - \sum_{i=0}^r 
\frac{1}{\beta_{ii}} \Big(b J_p(S_{p}^{i, k+1},\vec S_{\vec u}^{i, k}) - 
 J_d(\vec S_{\vec 
u}^{i, k+1},\vec S_{\vec u}^{i, k})\Big)\,.
\end{align}
For the second term $T_2$ on the right-hand side, \eqref{proof_eq_19} 
and \eqref{Eq:BoundJ} yield that
\begin{equation*}
|T_2| \le  c \sum_{i=0}^r 
\dfrac{b^2}{ \lambda \beta_{ii}} \| S_p^{i, k+1} - S_p^{i, k} \|^2 \,.
\end{equation*}
On the left-hand side, the third term can be rewritten by \cite[p.\ 760, Eq.\ 
(4.29)]{BRK17}. The fourth and fifth term are rewritten by first using the identity $2 
\langle x,y\rangle = \langle x,x\rangle + \langle y,y\rangle- \langle x-y,x-y\rangle$ and 
then applying \eqref{proof_eq_19} with \eqref{Eq:BoundJ}.

\smallskip
\noindent\textbf{7.\ Step (Contraction).} 
We combine the results of the 5th and 6th step to
\begin{align}
\nonumber 
& \sum_{i=0}^r \Big\{\Big(\dfrac{1}{M \beta_{ii}} + \dfrac{L}{2 \beta_{ii}}\Big) \| S_p^{i, k+1} \|^2  + \dfrac{L}{ 2 \beta_{ii}} \| S_p^{i, k+1} - S_p^{i, k} \|^2   \Big\}\\[0ex]
\nonumber
&  \; \quad +  \dfrac{\tau_n}{2} \|\vec K^{-1/2}  \vec e_{\vec q}^{k+1} (t_{n}^-)\|^2 + \dfrac{\tau_n}{2} \|\vec K^{-1/2}  \vec e_{\vec q}^{k+1} (t_{n-1}^+)\|^2\\[0ex]
\nonumber
&  \; \quad + \sum_{i=0}^r \frac{1}{2\beta_{ii}}\sum_{K\in \mathcal T_h} \la \vec \sigma (\vec S_{\vec u}^{i, k+1}), \vec \veps(\vec S_{\vec u}^{i, k+1})\ra_K + \sum_{i=0}^r \frac{1}{2\beta_{ii}} J_\delta (\vec S_{\vec u}^{i,k+1},\vec  S_{\vec u}^{i, k+1})\\
\nonumber 
&  \; \quad + \sum_{i=0}^r \frac{1}{2\beta_{ii}}\sum_{K\in \mathcal T_h} \la \vec \sigma (\vec S_{\vec u}^{i, k}), \vec \veps(\vec S_{\vec u}^{i, k})\ra_K +  \sum_{i=0}^r \frac{1}{2\beta_{ii}} J_\delta (\vec S_{\vec u}^{i,k},\vec  S_{\vec u}^{i, k})\\
\label{proof_eq_30}
& \leq  \; \sum_{i=0}^r  \dfrac{L}{ 2 \beta_{ii}} \| S_p^{i, k} \|^2 + c \sum_{i=0}^r 
\dfrac{b^2}{\lambda \beta_{ii}} \| S_p^{i, k+1} - S_p^{i, k} \|^2 \\
\nonumber 
&  \; \quad -  \sum_{i=0}^r \frac{b}{\beta_{ii}} 
J_p(S_{p}^{i, k+1},\vec S_{\vec u}^{i, k})
+ \sum_{i=0}^r  \frac{1}{\beta_{ii}}  J_d(\vec S_{\vec u}^{i, k+1} ,\vec S_{\vec u}^{i, k}) \,.
\end{align}
Using \eqref{Eq:BoundJ} the last two terms on the right-hand side of \eqref{proof_eq_30} 
can be absorbed by terms on the left-hand side, if $L>0$ and the penalty function $\delta$ 
in $J_\delta$ of \eqref{Eq:D_AlgebProb3} are chosen sufficiently large, i.e.\  $L > 
c b^2/ \lambda$.  Inequality \eqref{proof_eq_30} then shows the convergence of the 
iterates $S_p^{i, k}$ in 
$W_h$, for $i=0,\ldots, r$.  From \eqref{proof_eq_30} along with the convergence of $S_p^{i, k}$ we get the 
convergence of $\vec e_{\vec q}^{k} (t_{n}^-)$ and $\vec e_{\vec q}^{k} (t_{n-1}^+)$ to 
$\vec{0}$ for $k\rightarrow \infty$. Eq.\ \eqref{proof_eq_8} together with 
the convergence of $\vec e_{\vec q}^{k} (t_{n}^-)$ yields the convergence of 
$e_p^{k} (t_{n}^-)$ to $0$ for $k \rightarrow \infty$. Finally, Eq.\ 
\eqref{proof_eq_9} implies the convergence of $\vec e_{\vec u}^{k} (t_{n}^-)$ to 
$\vec {0}$. The rest follows as in \cite[Cor.\ 4.6]{BRK17}.
\end{mproof}

\vspace*{-1ex}

\ifx\undefined\bysame\newcommand{\bysame}{\leavevmode\hbox to3em{\hrulefill}\,}
\fi

\end{document}